\documentclass[conference]{IEEEtran}
\IEEEoverridecommandlockouts
\usepackage{cite}
\usepackage{amsmath, amssymb, amsfonts}
\usepackage{amsthm}
\usepackage{algorithmic}
\usepackage{graphicx, float,wrapfig,subfigure} 
\usepackage{svg} 
\usepackage{overpic} 
\usepackage{textcomp}
\usepackage{multirow}
\usepackage{xcolor}
\usepackage{nomencl}
\usepackage{ifthen} 
\usepackage[super]{nth} 
\usepackage{soul} 
\usepackage{url} 
\usepackage{footnote}
\allowdisplaybreaks

\newcommand{\red}[1]{{\color{red}#1}}
\newcommand{\blue}[1]{{\color{blue}#1}}

 \def\BibTeX{{\rm B\kern-.05em{\sc i\kern-.025em b}\kern-.08em
    T\kern-.1667em\lower.7ex\hbox{E}\kern-.125emX}}
    
\newcounter{remark}
\newenvironment{remark}[1][]{\refstepcounter{remark}\par
	\textit{Remark~\theremark #1:} \rmfamily}

\newcounter{proposition}
\theoremstyle{definition} 
\newtheorem{definition}{Definition}[section]

\begin{document}
\title{Solving Overstay in PEV Charging Station Planning via Chance Constrained Optimization}
\author{Teng Zeng, \IEEEmembership{Student Member, IEEE}, Hongcai Zhang, \IEEEmembership{Member, IEEE}, and Scott Moura, \IEEEmembership{Member, IEEE},\\
\vspace{-10 mm}
\thanks{This work was supported by the National Science Foundation through grant 1746858. T. Zeng, H. Zhang, and S. Moura are with the Department of Civil and Environmental Engineering, University of California, Berkeley 94720 CA (email:zhang-hc13@berkeley.edu).  S. J. Moura is also with the Smart Grid and Renewable Energy Laboratory, Tsinghua-Berkeley Shenzhen Institute, Shenzhen, 518055, P. R. China.}}
\maketitle

\begin{abstract}
     Currently, the utilization rate of public plug-in electric vehicle (PEV) charging stations is only 15\%. During the other 85\% of the time, the chargers are either idle or occupied by a fully charged PEV. Strategically switching a fully charged PEV to another one waiting to charge, which we refer to as \textit{interchange}, is an indispensable measure to enhancing charging station utilization. This paper studies the PEV charging station planning problem considering strategic interchange, which has not yet been well studied in the literature. The objective is to enhance utilization rate while lowering  capital investment and operation cost. A new power/energy aggregation model is proposed for interchange operation and a chance-constrained planning model with interchange is developed for a public charging station with customer demand uncertainties. 
     Numerical experiments are conducted to illustrate the performance of the proposed method. Simulation  results show that incorporating  strategic interchange operation, can significantly decrease the number of chargers, enhance utilization and economic  efficiency. 
\end{abstract}


\renewcommand{\nomgroup}[1]{%
\ifthenelse{\equal{#1}{D}}{\item[\textbf{Definitions/Abbreviations}]}{%
\ifthenelse{\equal{#1}{I}}{\item[\textbf{Indices/Sets}]}{%
\ifthenelse{\equal{#1}{V}}{\item[\textbf{Variables}]}{%
\ifthenelse{\equal{#1}{P}}{\item[\textbf{Parameters}]}{}}}}}

\small
\makenomenclature
        \nomenclature[D]{PEV}{Plug-in electric vehicle}
        \nomenclature[D]{SoC}{Battery state of charge}
        
        \nomenclature[I]{$i/ I$}{index of PEVs}
        \nomenclature[I]{$t,\tau/ T$}{index/set of specified time stamp of a day}
        
        \nomenclature[P]{$B_i$}{$i$\textsuperscript{th} PEV's battery capcity}
        \nomenclature[P]{$c^{\text{ch}}$}{ Per-unit cost for charger (\$) }
		\nomenclature[P]{$c_t^{\text{e}}$}{ Per-unit cost for electricity consumption (\$/kWh) at time $t$ }
		\nomenclature[P]{$c^{\text{loss}}$}{ Per-unit cost for load shedding (\$/kWh) }
		\nomenclature[P]{$c_{\text{plan}}^{\text{itc}}$}{Per-unit cost for plug-in energy demand interchange (\$/kWh) }
	    \nomenclature[P]{$c_{\text{oper}}^{\text{itc}}$}{Per-unit cost for interchange (\$) }
		\nomenclature[P]{$c^{\text{ed}}$}{Per-unit cost for electricity peak demand charge (\$/kWh) }
		\nomenclature[P]{$e_i^{\text{need}}$}{$i$\textsuperscript{th} PEV's required energy when departs } 
		\nomenclature[P]{$e_i^{\text{pneed}}$}{$i$\textsuperscript{th} PEV's minimum (required) plug-in consumption }
		\nomenclature[P]{$e_i^{\text{max}}$}{$i$\textsuperscript{th} PEV's maximum energy possibly acquired by its battery }
		\nomenclature[P]{$e_i^{\text{pmax}}$}{$i$\textsuperscript{th} PEV's maximum plug-in energy }
		\nomenclature[P]{$e_{i,\tau}^{\text{p}+/-}$}{$i$\textsuperscript{th} PEV's plug in energy upper/lower boundary at time $\tau$ } 
		\nomenclature[P]{$\eta$}{Efficiency during charging}
		\nomenclature[P]{$p_{i,\tau}^{\text{p}+/-}$}{$i$\textsuperscript{th} PEV's plug in power upper/lower boundary at time $\tau$ } 
		\nomenclature[P]{$p^{\text{rated}}$}{Charger's rated charging power } 
		\nomenclature[P]{$p_{\text{max/min}}^{\text{tran}}$}{Local transformer power upper/lower bound } 
		\nomenclature[P]{$SoC_i^{\text{d}}$}{$i$\textsuperscript{th} PEV's expected departure SoC}
		\nomenclature[P]{$t^\text{a/d}_i$}{$i$\textsuperscript{th} PEV's arrival/expected departure time}
        \nomenclature[P]{$\Delta t$}{Duration of each sub-hourly time interval}
		\nomenclature[P]{$\Delta t_i^{\text{itc}}$}{$i$\textsuperscript{th} PEV's interchange time}
		\nomenclature[P]{$\Tilde{P}_\text{max}^{\text{grid}}$}{Estimated grid peak power demand based on forecasting}
		\nomenclature[P]{$P_\text{max}^{\text{grid}*}$}{Past recorded net peak $P_\text{max}^{\text{grid}}$ of the month}
        
        
        \nomenclature[V]{$P_t^{\text{loss}}$}{PEV load shedding at time $t$}
        \nomenclature[V]{$P_t$}{Aggregate PEVs power profile at time $t$}
        \nomenclature[V]{$P_t^{\text{p}}$}{Aggregate PEVs plug-in power profile at time $t$}
        \nomenclature[V]{$P_t^{\text{pitc}}$}{Plug-in power demand shed by interchanges at time $t$}
		\nomenclature[V]{$P_t^{\text{grid}}$}{Net grid power at time $t$}
		\nomenclature[V]{$X$}{Number of chargers}
		\nomenclature[V]{$S_{i,t}$}{$i$\textsuperscript{th} PEV plug in status at time $t$, 1 for plugged in, 0 otherwise}
	    \nomenclature[V]{$S_{i,t}^{\text{itc}}$}{$i$\textsuperscript{th} PEV interchange status at time $t$, 1 for interchange, 0 otherwise}
		
\printnomenclature

\normalsize
\section{Introduction}
\IEEEPARstart{P}{lug-in} electric vehicles (PEVs) are more energy efficient than conventional internal combustion engine vehicles \cite{hardman2018review}\cite{poullikkas2015sustainable}. They are characterized by reduced greenhouse gas (carbon-dioxide), air pollutants and noise emissions. A recent study based on public charging station data forecasts that the anticipated number of PEVs will reach 1 million in the U.S market by 2020, and more than 50\% of new cars sold globally will be electrified by 2040 \cite{ChargePointReport}. However, the continuing growth of PEVs might be impeded by its limited driving range and lack of public charging infrastructure. Although governments and private parties have put forth great efforts to facilitate the growth of the public charging systems, there remains a large gap between the current service capability and the expected PEV deployment. That is, PEV penetration has out-paced charging station deployment \cite{chinaSpeedupEV}\cite{charginglocations}. Furthermore, due to improper  planning of station sizes and sites, the charging facilities experience significantly imbalanced utilization rate. In urban areas, especially central business districts, the competition for charging resources is intense. After a charger is plugged in, it can be occupied (even if the PEV is not charging) until the driver returns (from work, shopping, dining, etc.). Nowadays, this kind of charger occupation takes 6-8 hours in a typical day. The limited number of available chargers may result in unbearable queuing and inconvenience that can significantly degrade the quality of service \cite{kong2016smart}. 
Therefore, properly planning and dealing with over-staying PEVs becomes an important yet still poorly understood issue.




    

In literature, reference \cite{lindgren2015bottlenecks} identifies that fully recharged PEVs continuing to physically occupy the charging slots creates a service bottleneck. To deal with this issue, \cite{biswas2016managing} introduces a penalty function. 
The penalty is activated once the actual charging session is finished, but the PEV is still occupying the charger. The trade off between the penalty price and the user acceptance probability is examined. Recently, Tesla has implemented a similar approach to address this exact issue. They impose an ``idle fee," which is a penalty cost they apply to users, dollar by minute, when the PEV is fully charged \cite{lambertTesla}. Besides penalizing the delayed user behaviors, \cite{lindgren2015bottlenecks} and \cite{zhang2017optimal} proposed the idea of single pole multiple cables. This allows multiple PEVs to simultaneously connect to one charging circuit, but only one PEV is charged at a time. Then, once the PEV is fully charged, the power output is switched to another PEV. 




While operational ideas exist to address low utilization, most of the charging station planning literature only focuses on monetary utilities (maximizing profits/social welfare). 
The focus is typically mitigating traffic congestion in a transportation system or shaving the peak load for the distribution power system. 
For example, in \cite{wang2013traffic}, a multi-objective PEV charging station planning model considering traffic constraints was proposed. Reference \cite{yao2014multi} studied coordinated planning of the integrated power distribution network and PEV charging systems. References \cite{wang2013traffic} and \cite{yao2014multi} both assumed PEV charging demands to be proportional to the traffic flow, modeled as a Poisson distribution. Reference \cite{zhang2016integrated} develops an aggregated planning framework for various types of charging facilities in an urban area. 
To the authors' best understanding, literature has not yet addressed either the PEV overstaying problem or the low utilization rate from a planning perspective.


In this paper, we study PEV charging station planning taking charger interchange into account. The interchange concept is precisely defined in Section II. We find that strategic interchange can significantly enhance station utilization rates.  
In addition, we introduce chance-constrained optimization to handle the parameter uncertainties, such as human driver behaviors. We take the perspective of a commercial business seeking to provide charging services, who has budget constraints. Therefore, the planning objective focuses on destination charging stations and lowering the cost of the entire investment. This paper advances the relevant literature through the following contributions:
    \begin{itemize}
        \item We address the PEV overstaying issue and model interchange during charging sessions -- a novel advancement to PEV charging station modeling.
        \item We propose an aggregate model to describe the station charging demand profile represented as power and energy constraints. This aggregate model, extended from \cite{zhang2017evaluation}, is governed by the PEV fleet of plug-in power and energy instead of their actual power and energy consumption. This improved model better resembles real-life operation.
    \end{itemize}

The interchange mechanism is defined in Section II and the aggregate demand model is introduced in Section III. Sections IV formulates the interchange planning model and Section V presents case studies to illustrate the proposed planning methodology. Results are summarized in Section VI.

\section{PEV Interchange Mechanism}

Our proposed idea includes the essence of coordinated charging, which takes PEV drivers' trip demands and real time local distribution system electricity supply into account and then provides charging services to the end users accordingly. It preserves the mutual benefits of peak shaving to the power grid as well as lower charging cost to the customers \cite{hardman2018review}. Moreover, our idea includes an extension to ``coordianted charging", which we refer to as ``interchange:"
\begin{definition}[\textit{Interchange}]
        One interchange (ITC) event is one occurrence that a charger which has been occupied by PEV A is switched (either automatically by machine or manually by human) to another waiting PEV B, before PEV A's departure time. Mathematically, when $t^{c.\text{start}}_B \leq t^{\text{d}}_A$ holds, one interchange event occurs.
\end{definition}

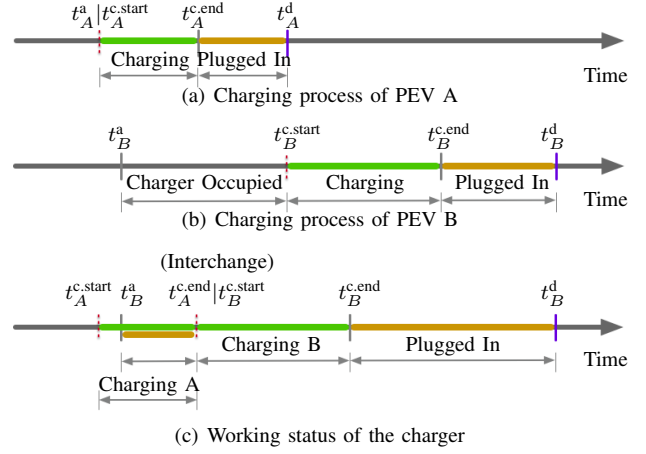
\begin{figure}
    \vspace{-4mm}
    	\centering
    	\subfigure[Charging process of PEV A]{
    		\begin{minipage}[h]{0.45\textwidth}
    			\begin{overpic}[width=1\columnwidth]{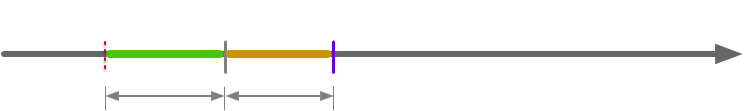}
    				\begin{normalsize}
    					\put(16,3.5){\footnotesize \color{black}\shortstack{Charging}}
    					\put(30,3.5){\footnotesize \color{black}\shortstack{Plugged In}}
    					\put(93,1.0){\footnotesize \color{black}\shortstack{Time}}
    					\put(9.5,11){\footnotesize \color{black} \shortstack{$t_{A}^{\text{a}}|t_{A}^{\text{c.start}}$}}
    					\put(27.5,11){\footnotesize \color{black}\shortstack{$t_{A}^{\text{c.end}}$}}
    					\put(43,11){\footnotesize \color{black}\shortstack{$t_{A}^{\text{d}}$}}
    				\end{normalsize}				
    			\end{overpic}
    		\end{minipage}
    	}
    	\subfigure[Charging process of PEV B]{
    		\begin{minipage}[h]{0.45\textwidth}
    			\begin{overpic}[width=1\columnwidth]{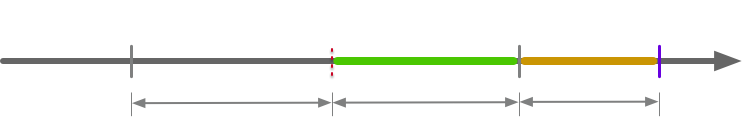}
    				\begin{normalsize}
    					\put(18.5,3.7){\footnotesize \color{black}\shortstack{Charger Occupied}}
    					\put(51,3.7){\footnotesize \color{black}\shortstack{Charging}}
    					\put(72,3.7){\footnotesize \color{black}\shortstack{Plugged In}}
    					\put(93,1.0){\footnotesize \color{black}\shortstack{Time}}
    					\put(15.5,11.5){\footnotesize \color{black}\shortstack{$t_{B}^{\text{a}}$}}
    				    \put(42.5,11.5){\footnotesize \color{black}\shortstack{$t_{B}^{\text{c.start}}$}}
    					\put(67.5,11.5){\footnotesize \color{black}\shortstack{$t_{B}^{\text{c.end}}$}}
    					\put(86,11.5){\footnotesize \color{black}\shortstack{$t_{B}^{\text{d}}$}}
    				\end{normalsize}				
    			\end{overpic}
    		\end{minipage}
    	}
    	\subfigure[Working status of the charger]{
    		\begin{minipage}[h]{0.45\textwidth}
    		\vspace{4mm}
    			\begin{overpic}[width=1\columnwidth]{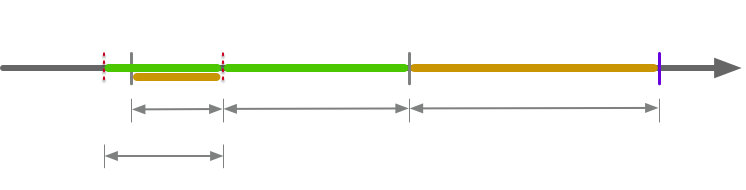}
    				\begin{normalsize}
    					\put(14,6){\footnotesize \color{black}\shortstack{Charging A}}
    					\put(34,12.5){\footnotesize \color{black}\shortstack{Charging B}}
    					\put(64,12.5){\footnotesize \color{black}\shortstack{Plugged In}}
    					\put(93,10){\footnotesize \color{black}\shortstack{Time}}
    					\put(8.5, 21){\footnotesize \color{black}\shortstack{$t_{A}^{\text{\text{c.start}}}$}}
    					\put(17.5,21){\footnotesize \color{black}\shortstack{$t_{B}^{\text{a}}$}}
    					\put(24,21){\footnotesize \color{black}\shortstack{(Interchange)\\$t_{A}^{\text{c.end}}|t_{B}^{\text{\text{c.start}}}$}}
    					\put(53,21){\footnotesize \color{black}\shortstack{$t_{B}^{\text{c.end}}$}}
    					\put(86,21){\footnotesize \color{black}\shortstack{$t_{B}^{\text{d}}$}}
    				\end{normalsize}				
    			\end{overpic}
    		\end{minipage}
    	}
    	\caption{A schematic diagram of two PEVs sharing one charger.}
    	\label{fig:sharing}
    	\vspace{-4mm}
\end{figure}

A schematic diagram of two PEVs A and B sharing one charger is shown in Figure \ref{fig:sharing}. Without proper management, PEV A occupies the charger after it is fully charged at $t_A^{\text{c.end}}$ and until PEV A departs at $t_A^{\text{d}}$. As a result, when PEV B arrives at $t_B^{\text{a}}$, it must wait until the charger is released. Note that the charger is block from providing charging service between $t_A^{\text{c.end}}$ and $t_A^{\text{d}}$. At a destination charging station, such as a workplace, the queued PEV may wait several hours. As a result, the infrastructure planner may deploy two chargers to satisfy the ongoing demands. However, if PEV B could be \textit{interchanged} at $t_A^{\text{c.end}}$ such that there is no or little idle time, then utilization rate and quality-of-service increases and we can avoid investing in additional chargers.

\section{PEV Aggregate Model}

This section introduces an aggregate model to describe a PEV fleet's aggregate plugged in power and energy demands in a charging station with interchange taken into account. We adopt and modify the model from \cite{zhang2017evaluation}, whose accuracy and efficient computational features are discussed within.

\subsection{Aggregate Model for PEV Charging Demand}
The PEV fleet's actual power and energy needs are given by the following constraints, respectively, 
\begin{align}
    & p_{t}^{-} \leq P_{t}\leq p_{t}^{+}, ~~\forall t,\label{eqn:power cumu}\\
	& e_{t}^{-}\leq \sum_{\tau=t_0}^{t}{P_{\tau}}\eta\Delta t\leq e_{t}^{+}, ~~\forall t,\label{eqn:energy cumu}
\end{align}
where $p_{t}^{-}$ and $p_{t}^{+}$ are the aggregate power lower and upper bounds based on ``delay as long as possible" charging and ``immediate" charging; likewise, $e_{t}^{-}$ and $e_{t}^{+}$ are the corresponding aggregate energy upper and lower bounds.

\subsection{Aggregate Model for PEV Plug-in Demand}
Considering PEV overstay, the number of chargers in a charging station is not determined by the PEVs' physical electricity consumption (charging demands), but determined by when the PEVs are plugged-in and occupy the charger (plug-in demands). To better address this fact, we define two new terms, i.e., \textit{plug-in power} and \textit{plug-in energy}, to distinguish conventional PEV charging power and energy consumption, as follows:
\begin{definition}[\textit{plug-in power}]
    A PEV's \textit{plug-in power} is the rated charging power of the plugged-in charger. It is the power that the PEV occupies but not necessarily consumes. Mathematically:
    \begin{align}
            & p_{i,t}^\text{p} = \begin{cases}
              p_i^\text{rated}, ~~\text{if PEV $i$ is plugged into a charger}, \\
            0, ~~\text{otherwise.}
            \end{cases} \label{eqn: plug-in power}
    \end{align}
\end{definition}

\begin{definition}[\textit{plug-in energy}]
    A PEV's \textit{plug-in energy} is the hypothetical energy consumed throughout the PEV's entire plug-in duration. Mathematically:
    \begin{align}
            &e_i^\text{p} = p_i^\text{rated}\eta\left(t_i^{d}-t_i^{a}\right). \label{eqn: plug-in energy}
    \end{align}
\end{definition}

Both the plug-in power and energy are combined to form a PEV's plug-in demands. Note that the plug-in power and energy upper-bound PEV $i$'s consumed power $P_{i,t}$ and energy $E_{i,t}$, respectively. Mathematically: $P_{i,t} \leq p_{i,t}^\text{p}$ and $E_{i,t} \leq e_i^\text{p}$ where $E_{t} \sum_{\tau=t_0}^{t}{P_{\tau}}\eta\Delta t$.

\begin{remark}
PEVs' plug-in demands are critical, as they are more relevant factors that affect a station's service capability rather than PEVs' actual power and energy consumption.
\end{remark}

In the following subsections, we propose an aggregate plug-in demand model for a fleet of PEVs. 
\subsubsection{Individual Plug-in Demand Model}
When a PEV $i$ is fully charged, it may take up to $\Delta t_i^{\text{itc}}$ time to interchange its charger to another PEV that has been waiting in the queue. Then, the PEV's minimum (required) plug-in energy need is:
\begin{align}
    & e_i^{\text{pneed}} = e_i^{\text{need}} + p^{\text{rated}}\eta\Delta t_i^{\text{itc}}, \label{eqn: plug-in energy min}
\end{align}
where the first term in (\ref{eqn: plug-in energy min}) is the plug-in energy required upon the PEV's departure time; the second term is the additional plug-in energy resulting from the interchange time delay, after the departing PEV is fully charged. Parameter $\Delta t_i^{\text{itc}}$ can be estimated from historical interchange data.\footnote{In scenarios where interchange is realized by automatic machinery, it could be as few as zero time whereas in scenarios where human interactions are involved, it may be as long as several minutes to several hours.}

The maximum plug-in energy of PEV $i$ is determined by its expected plug-in duration, which generally is the same period as the parking duration, and can be calculated as:
    \begin{align}
        & e_i^{\text{pmax}} = p^{\text{rated}}\eta\left(t_i^{\text{d}}-t_i^{\text{a}}\right).
    \end{align}
The plug-in power and energy boundaries are thus:
\begin{align}
    & e_{i, \tau}^{p+} = \begin{cases}
        e_i^{\text{pmax}}, ~~\tau > t_i^{\text{d}},\\
        \min\left(e_{i, \tau-1}^{p+} + p^{\text{rated}}\eta\Delta t, e_i^{\text{pmax}} \right), t_i^{\text{a}} < \tau \leq t_i^{\text{d}}\\
        0, ~~\tau \leq t_i^{\text{a}},
    \end{cases}\label{eqn: plug-in energy upper bound}\\
    & e_{i, \tau}^{p-} = \begin{cases}
        e_i^{\text{pneed}}, ~~\tau > t_i^{\text{d}},\\
        \max\left(0, e_i^{\text{pneed}} - p^{\text{rated}}\eta(t_i^\text{d}-\tau) \right), t_i^{\text{a}} < \tau \leq t_i^{\text{d}},\\
        0, ~~\tau \leq t_i^{\text{a}},
    \end{cases}\label{eqn: plug-in energy lower bound}\\
    & p_{i, \tau}^{p+} = \begin{cases}
        p^{\text{rated}}, ~~t_i^{\text{a}} < \tau \leq t_i^{\text{d}},\\
        0, ~~\tau > t_i^{\text{d}} \text{~or~} \tau \leq t_i^{\text{a}},
    \end{cases}\label{eqn: plug-in power upper bound}\\
    & p_{i, \tau}^{p-} = 0, ~~\forall\tau, \label{eqn: plug-in power lower bound}
\end{align}
\subsubsection{Aggregate Plug-in Demand Model} The aggregate plug-in model is also the summation of individual plug-in models over the same period. This model describes the feasible set of all possible aggregate PEV plug-in power and energy consumption trajectories. 


\begin{remark}
    Without interchange management, the PEVs will occupy the chargers until the driver comes back (from work, shopping, etc.). As a result, the PEVs' plug-in consumption will be equal to their maximum plug-in consumption, i.e., $P_{t}^{\text{p}} = \sum_i^I p_t^{\text{p}+}$ and $E_{t}^\text{p} = \sum_i^I e_t^{\text{p}+}$. 
    
\end{remark} 

\section{Planning Model}
In this section, we propose a chance-constrained PEV station planning model considering investment and operation costs. Specifically, we examine the balance between installing more chargers at the investment stage and adopting more interchange operations at the operation stage. The planning model aims to minimize the total economic costs of the charging station, which includes the initial capital cost for the chargers, and the operation costs for electricity, interchange, and load shedding (due to limited electricity supply capacity). Thus, the planning model is formulated as follows:
\begin{subequations}
    \begin{align}
        \begin{split}
            \min\quad &\zeta c^{\text{ch}} \blue{X} + 12c^{\text{ed}} \blue{P_{\text{max}}^{\text{grid}}} \\
            & + 365\sum_t \left(c_t^{\text{e}} \blue{P_t} + c^{\text{loss}} \blue{P_t^{\text{loss}}} + c_{\text{plan}}^{\text{itc}} \blue{P_t^{\text{pitc}}} \right) \Delta t
            \label{eqn: 4_obj}
        \end{split}\\
        \text{s.t.:}\quad& \red{p_t^-} \leq \blue{P_t} + \blue{P_t^{\text{loss}}} \leq \red{p_t^+}, ~~\forall t, \label{eqn: 4_1}\\
        & \red{e_t^-} \leq \sum_{\tau=t_0}^t \left( \blue{P_{\tau}} + \blue{P_{\tau}^{\text{loss}}} \right)\eta\Delta t \leq \red{e_t^+}, ~~\forall t, \label{eqn: 4_2}\\
        & 0 \leq \blue{P_t^{\text{pitc}}} \leq \red{p_t^{\text{p}+}}, ~~\forall t, \label{eqn: 4_3}\\
        & \red{e_t^{\text{p}-}} \leq \sum_{\tau=t_0}^t (\red{p_\tau^{\text{p}+}} - \blue{P_\tau^{\text{pitc}}}) \eta \Delta t \leq \red{e_t^{\text{p}+}}, ~~\forall t, \label{eqn: 4_4}\\
        & \blue{P_t} \leq \red{p_t^{\text{p}+}} - \blue{P_t^{\text{pitc}}} \leq p^{\text{rated}} \blue{X}, ~~\forall t, \label{eqn: 4_5}\\
        & p_{\text{min}}^{\text{tran}} \leq \blue{P_t} + \red{p_t^{\text{base}}} \leq p_{\text{max}}^{\text{tran}}, ~~\forall t, \label{eqn:4_6}\\
        & \blue{P_t} + \red{p_t^{\text{base}}} \leq \blue{P_{\text{max}}^{\text{grid}}}, ~~\forall t, \label{eqn:4_7}\\
        & 0 \leq \blue{X}, 0 \leq \blue{P_t}, 0 \leq \blue{P_t^{\text{loss}}}, 0 \leq \blue{P_t^{\text{pitc}}}, ~~\forall t, \label{eqn:4_8}
    \end{align}
\end{subequations}
where the optimization variables are colored blue and the stochastic parameters are colored red to enhance clarity and expose structure. Note that the decision variables $P_t^{\text{p}}$, $P_t^{\text{grid}}$ are defined as $P_t^{\text{p}} = \red{p_t^{\text{p}+}} - \blue{P_t^{\text{pitc}}}$ and $P_t^{\text{grid}} = \blue{P_t} + \red{p_t^{\text{base}}}$. After eliminating these two equality constraints, the optimization problem is formulated as shown above. The first term in the objective function (\ref{eqn: 4_obj}) is the annualized investment cost for the PEV chargers; the second term takes electricity demand charges into account (see e.g. \cite{PGEdemandcharge}); lastly, the set of three terms are the annual operation costs for electricity consumption, load shedding, and interchanges. Constraints (\ref{eqn: 4_1}) and (\ref{eqn: 4_2}) define the power and energy boundaries for the aggregate PEV charging profile, whereas constraints (\ref{eqn: 4_3}) and (\ref{eqn: 4_4}) define the plug-in power and energy boundaries in association with interchange. Constraint (\ref{eqn: 4_5}) determines the number of PEV chargers needed to satisfy the peak aggregate plug-in power demand; meanwhile it also ensures that the PEVs are only charged during their plug-in period. Constraint (\ref{eqn:4_6}) lower and upper bounds the total power demand from this charging station by the local transformer's capacity. Constraint (\ref{eqn:4_7}) determines the amount of demand charge and the last constraint defines the feasible domain sets of the variables.

In the above formulation, parameters $\red{p_t^+}$, $\red{p_t^-}$, $\red{e_t^+}$, $\red{e_t^+}$, $\red{p_t^{\text{p}+}}$, $\red{p_t^{\text{p}-}}$, $\red{e_t^{\text{p}+}}$, $\red{e_t^{\text{p}-}}$ and $\red{p_t^\text{base}}$ are stochastic. Each of these parameters can be associated with uncertain human behavior. We consider variance among different individuals, and thus variance exists in the aggregate demand modeling. Therefore, chance-constrained linear programming (CCLP) is applied to capture the uncertainties. All of the inequality constraints with stochastic parameters can be formulated as follows:
    \begin{align}
        \red{b} &\leq \textbf{a}^\intercal {\bf x},\\
        1-\epsilon &\leq \textbf{Pr}(\red{b} \leq \textbf{a}^\intercal {\bf x}), \label{cclp}
    \end{align}
where \red{$b$} is the uncertain parameter and $\epsilon$ is the reliability threshold. Since the parameter only shows up on the left hand side of the linear constraints, the above formulation (\ref{cclp}) can be reformulated into its deterministic counterpart if the cumulative distribution function (CDF) of the uncertain parameter \red{$b$} is known or well-estimated:
    \begin{align}
        \text{F}_\red{b}^{-1}(1-\epsilon) \leq \textbf{a}^\intercal {\bf x}, 
    \end{align}
where $\text{F}_\red{b}^{-1}$ is the inverse CDF of \red{$b$}. If parameter \red{$b$} follows Gaussian distribution, $\mathcal{N}(\bar{b}, \sigma^2)$, then the equivalent deterministic counterpart of the CCLP (\ref{cclp}), becomes:
    \begin{align}
        \bar{b} - \Phi^{-1}(\epsilon)\sigma \leq \textbf{a}^\intercal {\bf x},
    \end{align}
where $\Phi^{-1}(\epsilon)$ is the $\epsilon$-quantile of the standard normal distribution. Following this transformation, all of the above constraints (\ref{eqn: 4_1}) - (\ref{eqn:4_8}) can be reformulated into CCLPs.

\section{Case Study}
In this section, we perform Monte-Carlo simulation to evaluate the effectiveness of the proposed planning model, while considering coordinated charging and interchange. This is an online operation simulation, with the objective to optimally schedule the charging process for each PEV to minimize operation costs.

    \hl{\subsection{Case Overview and Parameter Settings}}
        Our case study involves a destination charging station, which serves 50 PEVs/day. We assume individual PEV demands follow normal distributions that are independent and identically distributed. Therefore, the aggregate demands also preserve properties of a normal distribution. The Monte-Carlo simulation is based on data from the National Household Travel Survey \cite{NHTS2001}, where $t_i^\text{a}$, $\text{SoC}_i^\text{a}$, $t_i^\text{d}$, and $\text{SoC}_i^\text{d}$ for each PEV $i$ are generated. The reliability threshold used is $\epsilon = 0.2$. We assume constant charging rate, $\text{P}^\text{rated}=$ 6.6 kW, and set $\text{SoC}_i^\text{max}=$ 0.95. We adopt a 24 kWh vehicle battery capacity and 0.14 kWh/km energy consumption rate from \cite{zhang2017optimal}. 
        
        The total investment cost for each charger, including hardware, computer controller (for coordinated charging), installation, human labor, etc., is assumed to be \$4000 \cite{RMIstationCosts}. 
        We consider a charger lifetime of 15 years and discount rate of 6\% \cite{schroeder2012economics}. The capital recovery factor is $\frac{r(1+r)^N}{(1+r)^N-1}$, in which a ratio is used to calculate the present value of an annuity. The total investment costs are thus converted to a stream of equal annual payments over $N$ years, which typically is the product's lifetime. The interchange cost $c_\text{oper}^\text{itc}$ is \$0.44 per interchange event and $c_\text{plan}^\text{itc}$ is thus \$0.0167/kWh. An industrial time-of-use (TOU) electricity tariff is adopted from PG\&E \cite{PGETOU}. 
        The penalty for unsatisfied charging demand ($c^\text{loss}$) is assumed to be five times the average electricity tariff price, \$1.2/kWh.

        Case 0 is the base case, where regular coordinated charging is applied without interchange. 
        Case 1 adds interchange at a cost $c_\text{oper}^\text{itc} =$ \$0.44/ITC, and Case 2 considers a very small cost $c_\text{oper}^\text{itc} =$ \$0.003/ITC, e.g. interchanging with a robot.
        
    \subsection{Planning Results and Analysis}
    
    \begin{savenotes}
    \begin{table}[ht] 
    \vspace{-3mm}
            \centering
            \caption{Planning Costs Saving Summary}
            \label{tab1}
            \begin{tabular}{|c|c|c|c|c|}
            \hline
                Scenarios & ITC Cost & Capital Costs & Piles & Cost Reduction\\
                \hline
                Case 0 & {\$0 (0\%)\footnote{\textit{($\cdot$) indicates the portion of the cost over the total investment process (OPEX + CAPEX).}}}  & {\$11119 (21.3\%)}  & 26 & 0\%\\
                \hline
                Case 1 & \$1846 (3.6\%) & \$7825 (15.4\%) & 19 & 3.0\%\\
                \hline
                Case 2 & \$60 (0.1\%) & \$2059 (4.7\%)& 5 & 16.6\%\\
                \hline
            \end{tabular}
            \label{table:planning summary}
            \vspace{-4mm}
    \end{table}
    \end{savenotes}

    \begin{figure}[ht] 
    \vspace{-2mm}
        \centering
        \includegraphics[scale=0.172]{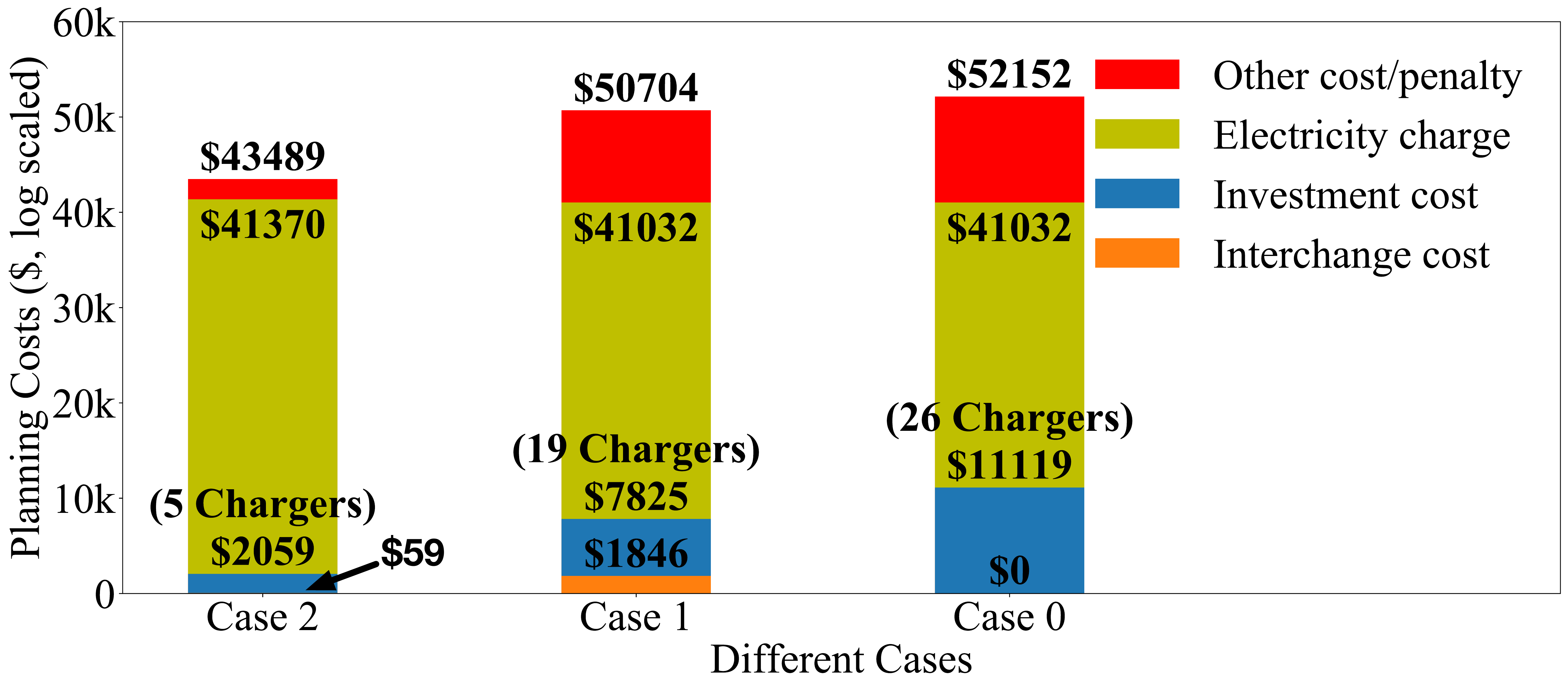}
        \vspace{-4mm}
        \caption{Planning Results Summary}
        
        \label{Planning summary1}
    \end{figure}
    
    The planning results are summarized in Table \ref{table:planning summary} and Fig. \ref{Planning summary1}. 
    Table \ref{table:planning summary} summarizes the major cost difference in the planning results. Adopting our proposed mechanism, more than one-third of the chargers are eliminated, resulting in a 3\% total cost savings (annual payment) in Case 1. Further improvements are shown in Case 2, where more than 80\% of the chargers are avoided and 16.6\% of the total cost is saved.
    
    \begin{figure}[htb]
    \vspace{-2mm}
        \centering
        \includegraphics[scale=0.11]{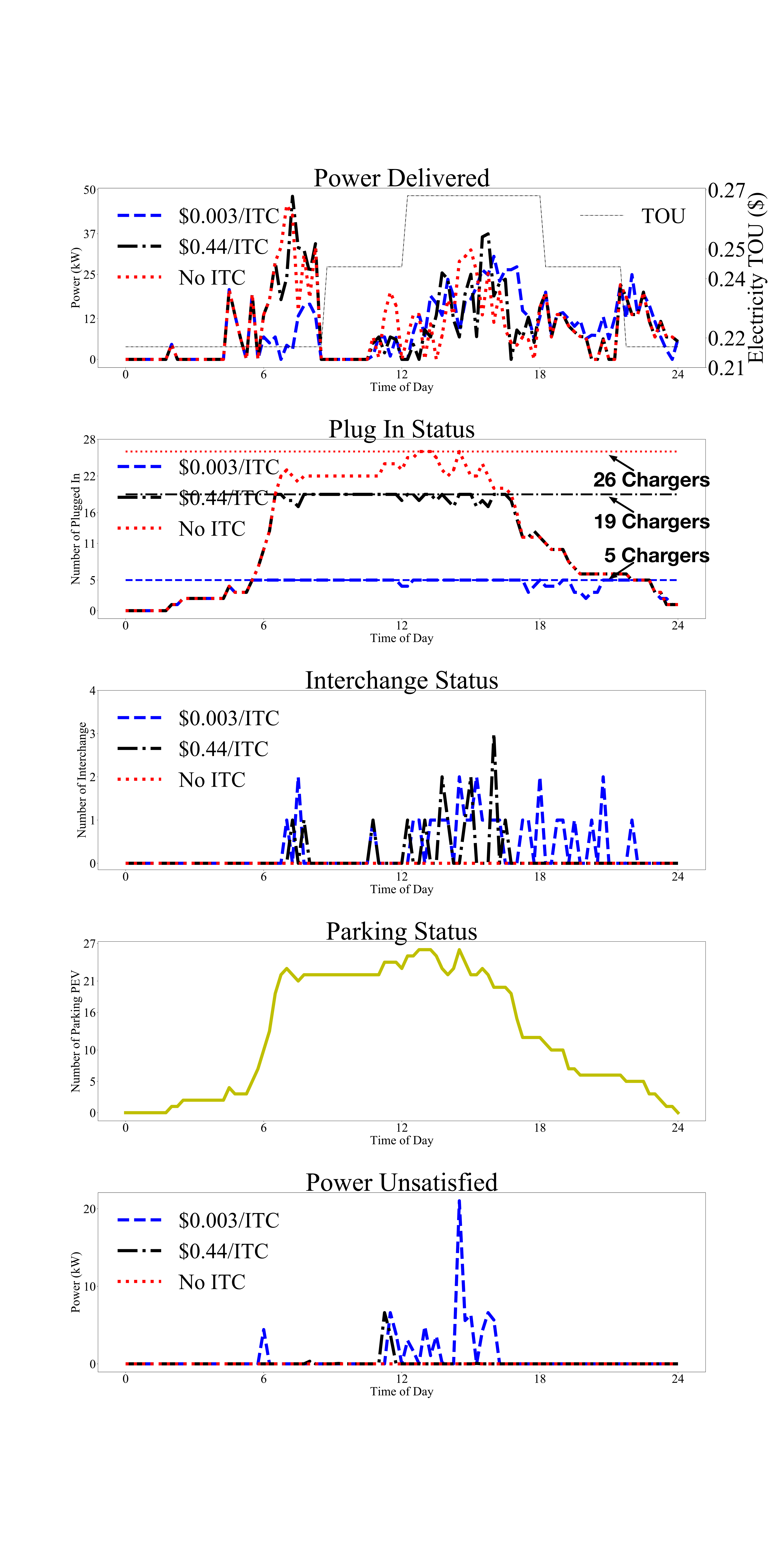}
        \vspace{-4mm}
        \caption{Rolling Horizon Operation}
        \label{Planning+operation1}
        \vspace{-4mm}
    \end{figure}
    Based on the planning results, one day of operation is simulated in Fig. \ref{Planning+operation1}. The simulation is rolling horizon-based. The station operator only optimizes the status of the chargers that are currently occupied by the PEVs present, meaning any future incoming PEV demands remain unknown at the current time. 
    During the peak demand period around 1pm, the station is visited by 26 PEVs.
    The entire aggregate arrival distribution throughout the day is plotted in the fourth subplot ``Parking Status" in Fig. \ref{Planning+operation1}. In Case 0, the station operator maxes out the number of installed chargers, 26 chargers, to avoid any unsatisfied demand. This is reflected in the second subplot of Fig. \ref{Planning+operation1}. In Case 2, only five chargers are needed due to cheap and frequent interchange events; the operation result from Fig. \ref{Planning+operation1} corroborates the feasibility and performance.

\section{Conclusion}



In this paper, we address optimal charging station planning in the face of the ``overstay'' problem.
We mathematically define and model vehicle ``interchange'', and accordingly propose an aggregate demand model. The planning results are illustrated through simulations. 
In one case study, we find that optimal interchange eliminates one-third of the chargers and reduces total cost by 3\%. If interchange events can be performed at nearly zero cost, then more than 80\% of chargers can be avoided and 16.6\% of costs are saved.
Extensions of this work can integrate other distributed energy resource assets (e.g. solar, storage), as well as bidirectional charging.

        




\bibliographystyle{unsrt}
\bibliography{bibliography.bib}

\end{document}